\documentclass[12pt]{article}

\usepackage{amsmath,amssymb,amscd}
\usepackage{cite}
\usepackage{graphicx,epstopdf}
\usepackage{float}
\usepackage{authblk} 
\numberwithin{equation}{section}
\usepackage[usenames]{color}
\usepackage{color}
\usepackage{colortbl}
\usepackage[linktocpage=true,plainpages=false,pdfpagelabels=false]{hyperref}
\definecolor{linkcolor}{rgb}{0.9,0,0}
\definecolor{citecolor}{rgb}{0,0.6,0}
\definecolor{urlcolor}{rgb}{0,0,1}
\hypersetup{
	colorlinks, linkcolor={linkcolor},
	citecolor={citecolor}, urlcolor={urlcolor}
}

\usepackage{soul}

\def\e{\epsilon}

\def\bp{\begin{proposition}}
\def\ep{\end{proposition}}
\def\bt{\begin{theo}}
\def\et{\end{theo}}
\def\be{\begin{equation}}
\def\ee{\end{equation}}
\def\bl{\begin{lemma}}
\def\el{\end{lemma}}
\def\bc{\begin{corollary}}
\def\ec{\end{corollary}}
\def\pr{\noindent{\bf Proof: }}

\def\bd{\begin{definition}}
\def\ed{\end{definition}}

\def\max{{\rm max\,}}

\def\max{{\rm max\,}}

\newtheorem{theo}{Theorem}[section]
\newtheorem{lemma}{Lemma}[section]
\newtheorem{definition}{Definition}[section]
\newtheorem{corollary}{Corollary}[section]
\newtheorem{proposition}{Proposition}[section]

\begin{document}

\title{Algebraic Geometry of Error Amplification: the Prony leaves}

\author[1]{Dmitry Batenkov}
\author[2]{Gil Goldman}
\author[3]{Yehonatan Salman}
\author[4]{Yosef Yomdin}

\affil[1] {Department of Mathematics, Massachusetts Institute of Technology, Cambridge, MA 02139, USA}
\affil[2,3,4]{Department of Mathematics, The Weizmann Institute of Science, Rehovot 76100, Israel}

\medskip

\affil[1]{Email: batenkov@mit.edu}
\affil[2]{Email: gilgoldm@gmail.com}
\affil[3]{Email: salman.yehonatan@gmail.com}
\affil[4]{Email: yosef.yomdin@weizmann.ac.il}

\maketitle

\begin{abstract}

We provide an overview of some results of \cite{Aki.Bat.Yom,Aki.Gol.Yom,Aki.Gol.Gol.Yom,Aki.Bat.Gol.Yom} on the
``geometry of error amplification'' in solving Prony system, in situations where the nodes near-collide.
It turns out to be governed by the ``Prony foliations'' $S_q$, whose leaves are ``equi-moment surfaces''
in the parameter space. Next, we prove some new results concerning explicit parametrization of the Prony leaves.

\end{abstract}


\section{Introduction}\label{Sec:Intro}
\setcounter{equation}{0}

We consider the problem of the ``measurements error amplification'' in solving classical
{\it Prony system} of algebraic equations, with the unknowns $a_j,x_j, \ j=1,\ldots,d,$ and with
the right hand side formed by the known ``noisy'' measurements $\mu_0,\ldots,\mu_{2d-1}$. This system has a form

\be\label{eq:Prony.system1}
\sum_{j=1}^d a_j x_j^k = \mu_k, \ k= 0,1,\ldots,2d-1.
\ee
We denote by $A=(a_1,\ldots,a_d) \in {\mathbb R}^d$ and $X=(x_1,\ldots,x_d) \in {\mathbb R}^d$
the unknowns in system (\ref{eq:Prony.system1}), and denote by ${\cal P}_d$ the ``parameter space'' of the unknowns $(A,X)$.
We will always assume that the nodes $X$ are pairwise different and ordered: $x_1<x_2<\ldots<x_d$.

\smallskip

Prony system appears in many theoretical and applied mathematical problems.
There exists a vast literature on Prony and similar systems - see, as a very small sample,
\cite{Pro,Aut,Bat1,Bat2,Bat.Yom2,Bey.Mon3,Pet.Plo,Pet.Pot.Tas,Plo.Wis,Pot.Tas} and references therein.
In particular, the bibliography in \cite{Aut} contains more than 50 pages.

Some applications of Prony system are of major practical importance, and,
in case when some of the nodes $x_j$ nearly collide, it is well known to present major mathematical
difficulties, in particular, in the context of ``super-resolution problem''
(see \cite{Aki.Bat.Yom,Aki.Gol.Yom,Aki.Gol.Gol.Yom,Aki.Bat.Gol.Yom,Aza.Cas.Gam,Bat1,Can.Fer1,Can.Fer2,Dem.Ngu,Dem.Nee.Ngu,Don,Fer,Mor.Can,Pet.Plo}
as a small sample).

\smallskip

The present paper deals with the problem of ``error amplification'' in solving a Prony system in the case that the nodes
$x_1,\ldots,x_d$ nearly collide. Our approach is independent of a specific method of inversion and deals with
a possible amplification of the measurements errors, in the reconstruction process,
caused by the geometric nature of the Prony system.

\smallskip

The paper consists of two parts: the first (Section \ref{Sec:Sum.Exs.Results}) is a summary of some recent results
of \cite{Aki.Bat.Yom,Aki.Gol.Yom,Aki.Gol.Gol.Yom,Aki.Bat.Gol.Yom} on the error amplification for near-colliding nodes.
The main observation here is that the incorrect reconstructions, caused by the measurements noise,
{\it are spread along certain algebraic subvarieties $S_q$ in the parameter space, which we call the ``Prony leaves''}.

\bd\label{def:Prony.leaf}
For $q=0,\ldots, 2d-1,$ and $\mu=(\mu_0,\ldots,\mu_q)$, the Prony leaf $S_q=S_q(\mu)$
is an algebraic variety in the parameter space ${\cal P}_d$,
defined by the first $q+1$ equations of the Prony system (\ref{eq:Prony.system1}):
\be\label{eq:Prony.system22}
\sum_{j=1}^d a_j x_j^k = \mu_k, \ k= 0,1,\ldots,q.
\ee
\ed
Generically, the dimension of the leaf $S_q(\mu)$ is $2d-q-1$. The chain
$$
S_0\supset S_1\supset \ldots \supset S_{2d-2}\supset S_{2d-1}
$$
can be explicitly computed (in principle), from the known measurements $\mu=(\mu_0,\ldots,\mu_{2d-1})$.
Notice that $S_{2d-1}$ coincides with the set of solutions of the ``full'' Prony system (\ref{eq:Prony.system1}).

\smallskip

In our approach the Prony leaves $S_q$ serve as an approximation to the set of possible ``noisy solutions'' of
(\ref{eq:Prony.system1}) which appear for a noisy right-hand side $\mu$.
The Prony curve $S_{2d-2}$ is especially prominent in the presentation below.

\smallskip

An important fact, found in \cite{Aki.Bat.Gol.Yom}, is that
{\it if the nodes $x_1,\ldots,x_d$ form a cluster of a size $h\ll 1$, while the measurements error is of order $\e$,
then the worst case error in reconstruction of $S_q$ is of order $\e h^{-q}$.
Thus, for smaller $q$, the leaves $S_q$ become bigger,
but the accuracy of their reconstruction becomes better.
The same is true for the accuracy with which $S_q$ approximate noisy solutions of (\ref{eq:Prony.system1}).}
Compare Theorems \ref{thm:main16} and \ref{thm:distance.to.Sq} below.

\smallskip

In particular, the worst case error in reconstruction of the solution $S_{2d-1}$ of (\ref{eq:Prony.system1}) is
$\sim \e h^{-2d+1},$ while the worst case error in reconstruction of the Prony curve $S_{2d-2}$
is of order $\e h^{-2d+2}.$ That is, the reconstruction of the Prony curve $S_{2d-2}$ is $h$ times
better than the reconstruction of the solutions themselves.

\smallskip

Consequently, we can split the solution of (\ref{eq:Prony.system1}) into two steps:
first finding, with an improved accuracy, the Prony curve $S_{2d-2}(\mu)$, and then localizing
{\it on this curve} the solution of (\ref{eq:Prony.system1}).
In particular, in the presence of a certain additional a priori information
on the expected solutions of the Prony system (for example, upper and/or lower bounds on the amplitudes),
it was shown in \cite{Aki.Bat.Gol.Yom} that the Prony curves
can be used in order to significantly improve the overall reconstruction accuracy.

\medskip

We believe that the results of \cite{Aki.Bat.Yom,Aki.Gol.Yom,Aki.Gol.Gol.Yom,Aki.Bat.Gol.Yom} presented in Section
\ref{Sec:Sum.Exs.Results} justify a detailed algebraic-geometric study of the Prony leaves.
In the second part of the present paper (Section \ref{Sec:Explicit.Prony.Param})
we prove some new results providing explicit equations, and explicit parametric representation, of the Prony leaves and their
projections into the nodes space.


\smallskip

Finally, in Section \ref{Sec:Open.Quest}, we summarize some open questions,
naturally arising in the study of the error amplification and of the Prony leaves.

\section{Setting of the problem}\label{Sec:Setting}
\setcounter{equation}{0}

In this paper we adopt one of many equivalent settings for the problem of inversion of the Prony system.
It is the problem of moment reconstruction of spike-trains, that is, of linear combinations of $d$ shifted
$\delta$-functions:
\be \label{eq:equation.model.delta}
F(x)=\sum_{j=1}^{d}a_{j}\delta\left(x-x_{j}\right),
\ee
with $A=(a_1,\ldots,a_d) \in {\mathbb R}^d, \ X=(x_1,\ldots,x_d) \in {\mathbb R}^d, x_1<x_2<\ldots <x_d$.
We will consider signal (\ref{eq:equation.model.delta}) as the point $(A,X)$ in the parameter space ${\cal P}_d$ introduced above.

We assume that the form (\ref{eq:equation.model.delta}) of signals $F$ is a priori known,
but the specific parameters $(A,X)\in {\cal P}_d$ are unknown.
Our goal is to reconstruct $(A,X)$ from $2d$ moments $m_k(F)=\int_{-\infty}^\infty x^k F(x)dx, \ k=0,\ldots,2d-1$,
which are known with a possible error bounded by $\e>0$.

\smallskip

An immediate computation shows that the moments $m_k(F)$ are expressed through
the unknown parameters $(A,X)$ as $m_k(F)=\sum_{j=1}^d a_j x_j^k$.
Hence our reconstruction problem is equivalent to solving the Prony system (\ref{eq:Prony.system1}), with
$\mu_k=m_k(F).$

\smallskip

Let a signal $F(x)\in {\cal P}_d$ be fixed.
In order to describe the geometry of the error amplification in solving (\ref{eq:Prony.system1}) we define,
following \cite{Aki.Gol.Gol.Yom,Aki.Bat.Gol.Yom},
the  $\e$-error set $E_\e(F)\subset {\cal P}_d$.
It consists of all signals $F'(x)\in {\cal P}_d$ which may appear in the reconstructions of
$F$ from noisy moment measurements $\mu'_k$, $|\mu'_k - m_k(F)|\le \e, \ k=0,\ldots,2d-1.$
Formally we have the following definition:

\bd\label{def:error.set}
The error set $E_\e(F)\subset {\cal P}_d$ is the set consisting of all the signals $F'(x)\in {\cal P}_d$ with
\be\label{eq:error.set}
|m_k(F')-m_k(F)|\le \e, \ k=0,\ldots,2d-1.
\ee
\ed
Our ultimate goal is a detailed understanding of the geometry of the error set $E_\e(F)$,
in the cases where the nodes of $F$ near-collide, and applying this information in order to improve the reconstruction accuracy.

\smallskip

We can explicitly describe the $\e$-error set $E_\e(F)$, considering the moments $m_k=m_k(F'), \ k=0,\ldots,2d-1,$
as non-linear coordinates in the space ${\cal P}_d$ of signals $F'$. Indeed, inequalities (\ref{eq:error.set})
immediately show that {\it in these coordinates the error set $E_\e(F)$ is the coordinate $\e$-cube in ${\cal P}_d$, centered at $F$}.

\smallskip

However, there are serious difficulties with this description. First, at the points where the nodes collide,
the moment coordinates develop complicated singularities. In particular, they fail to form a coordinate system near these points.
Still, the description of the error set $E_\e(F)$ via algebraic inequalities (\ref{eq:error.set}) remains valid.

\smallskip

Secondly, in the case of the nodes $X$ forming a cluster of size $h\ll 1$,
the ``moment coordinate system'' turns out to be significantly ``stretched'' in some directions,
up to the order $(\frac{1}{h})^{2d-1}$. Therefore the description of $E_\e(F)$ in the ``moment coordinate system''
given by (\ref{eq:error.set}) requires a ``translation'' into the natural coordinates $(A,X)$ in the space ${\cal P}_d$.
This translation (in a certain neighborhood of $F$ in ${\cal P}_d$) is
the main result of \cite{Aki.Bat.Gol.Yom} (see a review in Section \ref{Sec:Sum.Exs.Results} below).

\smallskip

We can aim to understand {\it the global geometry} of $E_\e(F)$
via an algebraic-geometric investigation of the moment coordinates, and, in particular, of the Prony
leaves.
In the case of two nodes, this investigation was started in \cite{Aki.Gol.Gol.Yom,Aki.Bat.Gol.Yom}.
In Section \ref{Sec:Explicit.Prony.Param} below we extend the results of \cite{Aki.Gol.Gol.Yom} to a general case of $d$ nodes.

\section{Summary of the results in \cite{Aki.Bat.Yom,Aki.Gol.Yom,Aki.Gol.Gol.Yom,Aki.Bat.Gol.Yom}}\label{Sec:Sum.Exs.Results}
\setcounter{equation}{0}

\medskip

Let a signal $F\in {\cal P}_d$ be given. We denote $I_F=[x_1,x_d]$ the minimal interval in
$\mathbb R$ containing all the nodes $x_1,\ldots,x_d$. We put $h(F)=\frac{1}{2}(x_d-x_1)$ to be the half of the length of $I_F$,
and put $\kappa(F)=\frac{1}{2}(x_1+x_d)$ to be the central point of $I_F$.

\smallskip

We ``normalize'' the signal $F$, shifting the interval $I_F$ to have its center at the origin,
and then rescaling $I_F$ to the interval $[-1,1]$. For this purpose we consider, for each $\kappa \in {\mathbb R}$ and $h>0$ the transformation

\be\label{eq:Ptransform}
\Psi_{\kappa,h}:{\cal P}_d\to {\cal P}_d,
\ee
defined by $(A,X)\to (A,\tilde X),$ with
$$
\tilde X=(\tilde x_1,\ldots,\tilde x_d), \ \ \tilde x_j=\frac{1}{h}\left (x_j-\kappa\right), \ j=1,\ldots,d.
$$
For a given signal $F$ we put $h=h(F), \ \kappa=\kappa(F)$ and call the signal $G=\Psi_{\kappa,h}(F)$ the model signal for $F$.
Clearly, $h(G)=1$ and $\kappa(G)=0$. Explicitly $G$ is written as
$$
G(x)=\sum_{j=1}^{d}a_{j}\delta\left(x-\tilde x_{j}\right).
$$
With a certain misuse of notations, we will denote the space ${\cal P}_d$ containing the model signals $G$ by $\tilde {\cal P}_d,$ and call it ``the model space''.
For $F\in {\cal P}_d$ and $G=\Psi_{\kappa,h}(F)$, the moments of $G$

\be\label{model.moments}
\tilde m_k(F)= m_k(G)=\sum_{j=1}^d a_j\tilde x_j^k, \ k=0,1,\ldots
\ee
are called the model moments of $F$.

\smallskip

\noindent{\it Below we describe the error set of $F$ in the associated model space $\tilde {\cal P}_d,$ and use the associated model moments}.

\smallskip

The main reason for mapping a general signal $F$ into the model space is that in the case of nodes $X$,
forming a cluster of size $h\ll 1$, as it was mentioned above, the moment coordinates turn out to be ``stretched'' in some directions,
up to the order $(\frac{1}{h})^{2d-1}$. {\it In contrast, in the model space $\tilde {\cal P}_d$ the system $m_0,\ldots, m_{2d-1}$
is compatible with the standard coordinates $(A,\tilde X)$ of $\tilde {\cal P}_d$, for all signals $G$ with
``well-separated nodes'' (see Theorem \ref{thm:coord.moments} below).}

\smallskip

For a given $F\in {\cal P}_d$, with the model signal $G=\Psi_{\kappa,h}(F),$
we denote by $\tilde E_\e(F)$ the set $\Psi_{\kappa,h}(E_\e(F)),$
which represents the error set $E_\e(F)$ of $F$ in the model space $\tilde {\cal P}_d$.
Note that $\tilde E_\e(F)$ is simply a translated and rescaled version of $E_\e(F)$
\smallskip

For given $\e,h>0$ denote by $\Pi_{\e,h}(G)$ the ``curvilinear parallelepiped'' consisting of all
$G' \in \tilde{\cal P}_d$ satisfying
the inequalities $$ |m_k(G')-m_k(G)|\leq \e h^{-k}, \ k=0,\ldots,2d-1.
$$

\bt\label{thm:main15}
For any $F\in {\cal P}_d$, let $\kappa=\kappa(F)$ and  $h=h(F)$. Let $G=\Psi_{\kappa,h}(F)$ be the model signal for
$F$.
Then, for any $\e>0$ we have $$
\Pi_{\e'',h}(G)\subset \tilde E_\e(F) \subset \Pi_{\e',h}(G),
$$
where $\e'=(1+|\kappa|)^{2d-1}\e, \ \e''=(1+|\kappa|)^{-2d+1}\e$. Specifically, for $\kappa=\kappa(F)=0$,
$$
\tilde E_\e(F)=\Pi_{\e,h}(G).
$$
\et
The result of Theorem \ref{thm:main15} holds without any assumptions on the mutual relation of $\e$ and $h$,
or on the distances between the nodes of $F$. However, without such assumptions it is difficult to provide any specific
geometric information on the moment parallelepiped $\Pi_{\e,h}(G)$. Still, Theorem \ref{thm:main15} implies an important fact:
{\it the Prony leaves $S_q$ of the model signal $G$ globally form a ``skeleton'' of the error set $\tilde E_\e(F)$, and,
in case when $\e$ and $h$ tend to zero in a certain rate, $S_q$ are the limits of $\tilde E_\e(F)$}.

\smallskip

To formulate this result accurately, let us denote the Prony leaves $S_q$ passing through $G$ by $S_q(G)$.
Thus $S_q(G)=S_q(\mu)$ for $\mu=(\mu_0,\ldots,\mu_{q})$ with $\mu_k=m_k(G), \ k=0,\ldots,q.$

\smallskip

Now let us assume that a model signal
$$
G=\sum_{j=1}^{d}a_{j}\delta\left(x-\tilde x_{j}\right) \in \tilde {\cal P}_d
$$
is fixed. For each real $\kappa$ and $h>0$ we consider a signal $F_{\kappa,h}(G)$, obtained from $G$ by an $h$-scaling and $\kappa$-shift of $x$:
$$
F_{\kappa,h}(G)=\Psi^{-1}_{\kappa,h}(G)=\sum_{j=1}^{d}a_{j}\delta\left(x- x_{j}\right) \in {\cal P}_d, \ x_j=h\tilde x_j+\kappa.
$$
Thus $G=\Psi_{\kappa,h}(F_{\kappa,h}(G))$ remains the model signal for each $F_{\kappa,h}(G)$.

For each $0 \le q \le 2d-1$ and $c>0$ denote by $S_q(G,c)$ the part of the Prony leaf $S_q(G),$
consisting of $G'\in S_q(G)$ with $|m_{q+1}(G')-m_{q+1}(G)|\le c.$

\begin{figure}
		\centering
		\includegraphics[scale=0.65]{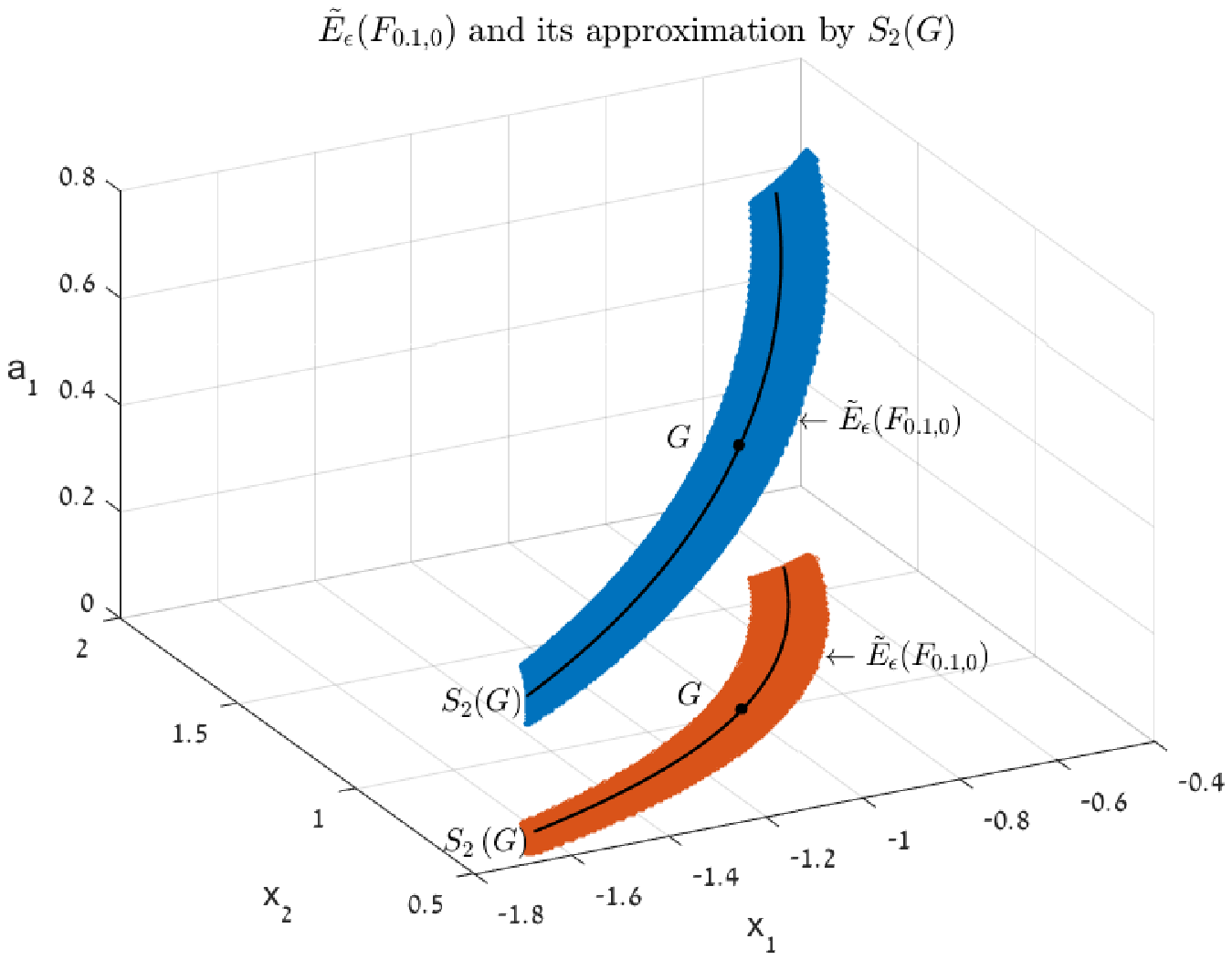}
		\caption{Presented are projections of the error set $\tilde{E}_{\e}(F_{h,\kappa})$ and a section of the Prony
		curve $S_2(G)$, for $G= \frac{1}{2} \delta\left(x+1 \right) + \frac{1}{2} \delta\left(x-1 \right)$,
		$h=0.1, \;  \kappa= 0$ and $\epsilon = h^3$. Stretched upwards is the
		projection into the coordinate subspace of $x_1,x_2,a_1$ and on the bottom plane into the nodes subspace
		$x_1,x_2$.}
		\label{fig.h01}
		\hspace{0.5cm}
		\centering
		\includegraphics[scale=0.65]{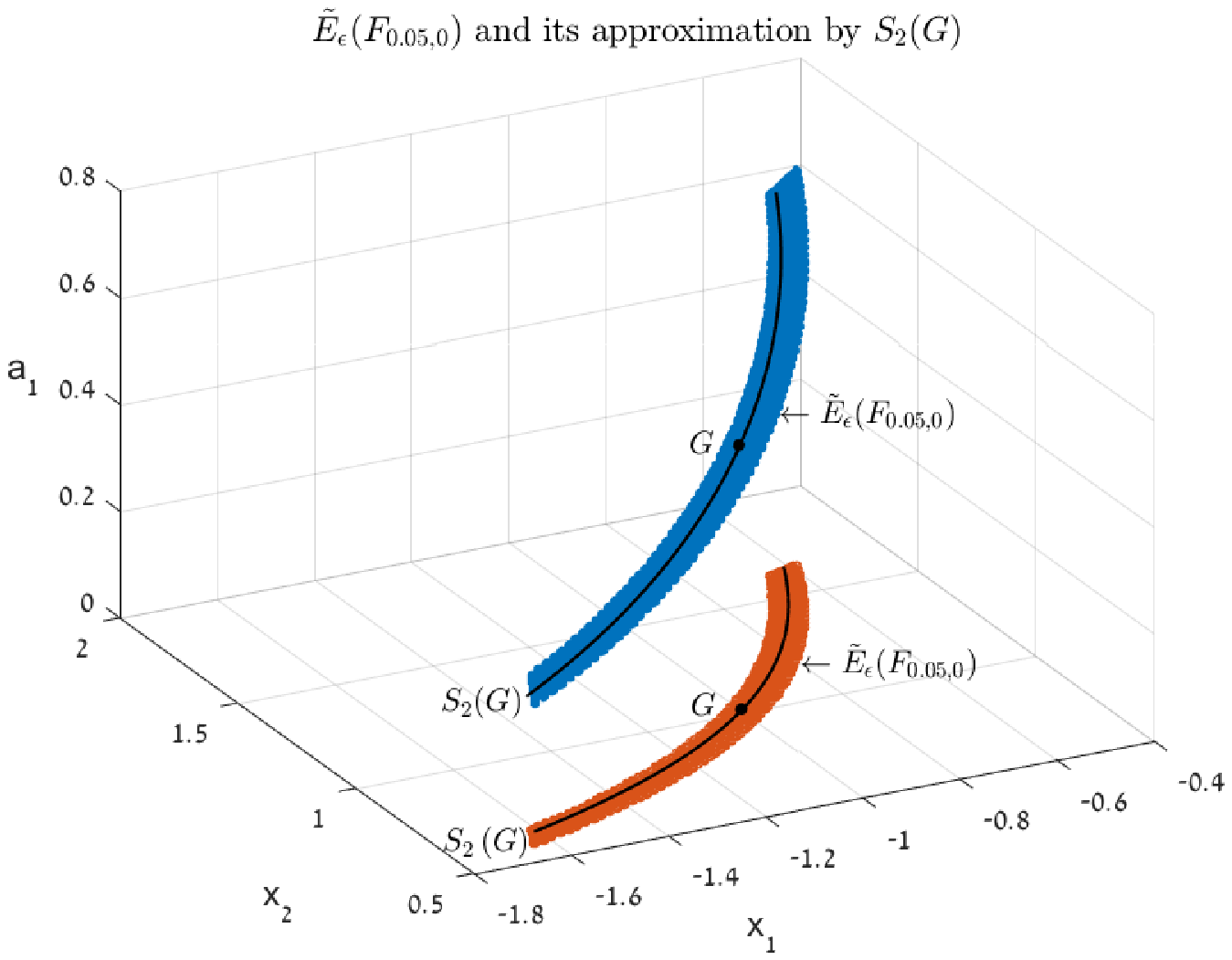}
		\caption{Presented are the error set $\tilde{E}_{\e}(F_{h,\kappa})$ and a section of $S_2(G)$ for
		 $G= \frac{1}{2} \delta\left(x+1 \right) + \frac{1}{2} \delta\left(x-1 \right)$, $h=0.05, \;  \kappa= 0$
		and $\epsilon = h^3$. Note the convergence of
		$\tilde{E}_{\e}(F_{h,\kappa})$ to $S_2(G)$.}
		\label{fig.h05}
\end{figure}

\bt\label{thm:Sq.limit}
Let $0 \le q \le 2d-1$, $\kappa$ and $C>0$ be fixed. Then for $h\to 0$ and for
$\e=Ch^{q+1},$ the error set $\tilde E_\e(F_{\kappa,h}(G))=\Psi_{\kappa,h}(E_\e(F_{\kappa,h}(G)))$ converges to the part $\bar S_q(G)$ of the Prony leaf $S_q(G),$
satisfying
$$
S_q(G,c'')\subset \bar S_q(G) \subset S_q(G,c'), \ \ c'=C(1+|\kappa|)^{2d-1}	, \ c''=C(1+|\kappa|)^{-2d+1}.
$$
\et
Figures \ref{fig.h01} and \ref{fig.h05} illustrate the case $d=2, q=2d-2=2$ of Theorem \ref{thm:Sq.limit}.

This theorem shows that the Prony leaves $S_q(G)$ globally approximate the error set
$\tilde{E}_\e(F)$, for $h\ll 1$ and $\e \sim h^{q+1}$. We consider the study of the ``before limit'' accuracy of this approximation as an important open question, which, presumably,
can be treated with the tools of real algebraic geometry. Some initial results in this direction, obtained in \cite{Aki.Bat.Gol.Yom},
and based on a ``quantitative'' version of the inverse function theorem, are presented below.

\smallskip

In order to apply this theorem, we have to make explicit assumptions on the separation of the nodes $X$ of our signal $F$, and on the size of its amplitudes $A$:

\bd\label{def:eta.regular}
Let $\eta$ satisfying $0<\eta\leq \frac{2}{d-1}$, $d>1$, and $m,M$ with $0<m<M$, be given.
A signal $G\in \tilde {\cal P}_d$ is called $(\eta,m,M)$-regular if for each $j=1,\ldots,d-1$ the distance between the neighbor nodes $\tilde x_j,\tilde x_{j+1}$ of $G$
is at least $\eta$, and the amplitudes $a_1,\ldots,a_d$ satisfy $m\leq |a_j|\leq M, \ j=1,\ldots,d$.

\smallskip

We say that a signal $F\in {\cal P}$ is $(h,\kappa,\eta,m,M)$-regular, if it can be obtained from an $(\eta,m,M)$-regular signal $G$ by an $h$-scaling,
and then a shift by $\kappa$.
\ed
We want to show that for an $(\eta,m,M)$-regular signal $G\in \tilde {\cal P}_d$ the model moments $m_0,\ldots,m_{2d-1}$ indeed form a coordinate system near $G$, which agrees with the standard coordinates $A,\tilde X$ on $\tilde {\cal P}_d$:

\bd\label{def:moment.coord.dist}
The moment metric $d(G',G'')$ on $\tilde {\cal P}_d$ is defined through the model moments $m_0,\ldots,m_{2d-1}$ as
$$
d(G',G'')=\max_{k=0}^{2d-1}|m_k(G'')-m_k(G')|.
$$
\ed

\bt\label{thm:coord.moments}
Let $G\in \tilde {\cal P}_d$ be an $(\eta,m,M)$ regular signal. Then there are constants $R,C_1,C_2$, depending only on $d,\eta,m,M,$ such that:

\smallskip

\noindent 1. The model moments $m_k=m_k(G')$ form a regular analytic coordinate system on the ball $B_R(G),$ centered at $G$,
of radius $R$ in the Euclidean metric on $\tilde {\cal P}_d$.

\smallskip

\noindent 2. The moment metric $d(G',G'')$ is Lipschitz equivalent on $B_R(G)$ to the Euclidean metric $||G''-G'||$: for each $G',G'' \in B_R(G)$ we have
$$
C_1 \ d(G',G'')\le ||G''-G'||\le C_2 \ d(G',G'').
$$
\et
Consequently, our description of the error set in terms of the model moments, given in Theorems \ref{thm:main15} and \ref{thm:Sq.limit} above, can be translated, inside the ball $B_R(G),$ into a description in the standard coordinates:

\bt\label{thm:distance.to.Sq}
Let $F$ be an $(h,\kappa,\eta,m,M)$-regular signal, and let $G\in \tilde {\cal P}_d$ be the model signal of $F$. Then
for each $q=0,\ldots,2d-1$ the ``local'' error set $\tilde E_\e(F)\cap B_R(G)$ is contained in the
$\Delta_q$-neighborhood (in the Euclidean metric) of the Prony leaf $S_q(G)$, for

$$
	\Delta_q=C_2\left(\frac{1+|\kappa|}{h}\right)^{q}\e.
$$

\et
Thus, as it was stated in the introduction, for smaller $q$, the leaves $S_q(G)$ become bigger, but the accuracy with
which they approximate the error set $\tilde E_\e(G)$ becomes better. The same is true for the accuracy of the
reconstruction of the leaves $S_q(G)$ from the noisy measurements: compare with Theorem \ref{thm:main16} below.

\subsection{Worst case reconstruction error}

We now present lower and upper bounds for the worst case reconstruction error $\rho(F,\e),$ defined by
$$
\rho(F,\e)= \max_{F'\in E_\e(F)}||F'-F||.
$$
In a similar way we define $\rho^A(F,\e)$ and $\rho^X(F,\e)$ - the worst case errors in reconstruction of the amplitudes $A=(a_1,\ldots,a_d)$ and of the nodes $X=(x_1,\dots,x_d)$ of $F$: considering signals $F'$ with parameters $A',X'$ we put
$$
\rho^A(F,\e)= \max_{F'\in E_\e(F)}||A'-A||, \ \rho^X(F,\e)= \max_{F'\in E_\e(F)}||X'-X||.
$$
\bt\label{thm:main6}
Let $F\in {\cal P}_d$ be an $(h,\kappa,\eta,m,M)$-regular signal. Then there are constants
$C_3,K_1,K_2,K_3,K_4$ depending only on $d,\kappa,\eta,m,M,$ such that for each positive $\e\le C_3h^{2d-1}$ the following bounds for the worst case reconstruction errors are valid:
$$
K_1 \e h^{-2d+1}\le \rho(F,\e), \ \ \rho^A(F,\e) \le K_2 \e h^{-2d+1},
$$
$$
K_3 \e h^{-2d+2}\le \rho^X(F,\e) \le K_4 \e h^{-2d+2}.
$$
\et
Theorem \ref{thm:main6} shows that the noise level $\e=\e_0=C_3h^{2d-1}$ plays a role of a threshold in noisy reconstruction: for $\e < \e_0$ the worst case error remain bounded and decreases with $\e$. However, for $\e \sim \e_0$ and bigger the nodes may collide and the amplitudes may blow up to infinity.

\smallskip





For $F\in {\cal P}_d$  and $G$ the model signal of $F$, we define the worst case reconstruction error of the Prony
leaves $S_q(G)$, in the model space ${\cal\tilde{P}}_d$, as follows.
For $q=0,\ldots,2d-1$
$$
	\rho^{S_q}(F,\e)= \max_{G'\in \tilde E_\e(F)} d_H(S_q(G)\cap B_R(G),S_q(G')\cap B_R(G)).
$$
Here $d_H(W_1,W_2)$ is the Hausdorff distance between the sets $W_1,W_2$, and $B_R(G)\subset \tilde {\cal
P}_d$ is the ball defined in Theorem \ref{thm:coord.moments}.

\bt\label{thm:main16}
Let $F\in {\cal P}_d$ be an $(h,\kappa,\eta,m,M)$-regular signal, and let $G$ be the model of $F$. Then, for each positive $\e\le C_3h^{2d-1}$ the following bounds for the worst case reconstruction errors of the Prony leaves $S_q(G), \ q=0,\ldots,2d-1$, are valid:
$$
K_1 \e h^{-q}\le \rho^{S_q}(F,\e) \le K_2 \e h^{-q}.
$$
\et



\section{Explicit parametrization of Prony leaves}\label{Sec:Explicit.Prony.Param}
\setcounter{equation}{0}

In this section we show that the Prony leaves $S_q(\mu)$ allow for an explicit parametrization.
For $q\le d-1$ (Section \ref{Sec:P.leaves.q.small}) this parametrization is produced in a rather straightforward way,
via expressing some of the amplitudes $a_j$ through the nodes and the remaining amplitudes.
It requires solving linear systems with Vandermonde matrices on the nodes $x_j$.
For $q\ge d$ (Section \ref{Sec:P.leaves.q.big}) our parametrization is produced via a proper modification of the classical solution method of the Prony system (suggested,
essentially, already in \cite{Pro}). It requires solving linear systems with Hankel matrices on the moments $\mu_k,$
and subsequently finding the roots of a univariate polynomial. In Section \ref{Sec:two.nodes} we illustrate the results of Sections
\ref{Sec:P.leaves.q.small} and \ref{Sec:P.leaves.q.big}, providing a complete description of the Prony leaves in the case of two nodes.

\subsection{Prony leaves $S_q$ for $q\leq d-1$}\label{Sec:P.leaves.q.small}

As above, we consider signals $F=(A,X)=\sum_{j=1}^d a_j\delta(x-x_j)\in {\cal P}_d$.
We denote by ${\cal P}^A_d$ and ${\cal P}^X_d$ the spaces of the amplitudes $A=(a_1,\ldots,a_d)$ and of the nodes $X=(x_1,\ldots,x_d),$
respectively, and denote by $\pi$ the projection
$$
\pi: {\cal P}_d\cong {\cal P}^A_d \times {\cal P}^X_d\to {\cal P}^X_d.
$$
For a given $\mu=(\mu_0,\ldots,\mu_q)$ we consider the Prony leaf $S_q(\mu)\subset {\cal P}_d,$

\bt\label{Thm:Prony.leaves.small.q}
For $q\leq d-1$ and for any $\mu=(\mu_0,\ldots,\mu_q)$ the Prony leaf $S_q(\mu)$ is a smooth subvariety of ${\cal P}_d$ of dimension $2d-q-1$.

\smallskip

The projection $\pi: S_q(\mu)\to {\cal P}^X_{d}$ is onto, and forms a regular locally trivial fibration over ${\cal P}^X_{d}$.
The fibers of $\pi$ are affine subvarieties of dimension $d-q-1$ in ${\cal P}_d$.

\smallskip

The amplitudes $a_{q+2},\ldots,a_d$ and all the nodes $x_1<x_2<\ldots <x_d$ can be chosen as the coordinates on $S_q(\mu)$,
while the amplitudes $a_1,\ldots,a_{q+1}$ are expressed in these coordinates as
$$
a_j=\sum_{l=0}^q A_l(X)\mu_l+\sum_{s=q+2}^d B_s(X)a_s, \ j=1,\ldots,q+1,
$$
with $A_l(X), B_s(X)$ regular rational functions in $X=(x_1,\ldots,x_d).$
\et
\pr


In the coordinates $(A=(a_1,\ldots,a_d),X=(x_1,\ldots,x_d))$ in ${\cal P}_d$ the Prony leaf $S_q(\mu)$ is defined by the following equations:

\begin{equation}\label{eq:Prony.fol.new.1}
\begin{array}{c}
a_1+a_2+\ldots + a_d=\mu_0\\
a_1x_1+a_2x_2+\ldots + a_dx_d=\mu_1\\
a_1x^2_1+a_2x^2_2+\ldots + a_dx^2_d=\mu_2\\
..........\\
a_1x^{q}_1+a_2x^{q}_2+\ldots + a_dx^{q}_d=\mu_{q}\\
\end{array}
\end{equation}


We can rewrite equations (\ref{eq:Prony.fol.new.1}) as

\begin{equation}\label{eq:Prony.fol.new.2}
\begin{array}{c}
a_1+a_2+\ldots+a_{q+1}=\mu_0-a_{q+2}-\ldots - a_d\\
a_1x_1+a_2x_2+\ldots+a_{q+1}x_{q+1}=\mu_1-a_{q+2}x_{q+2}-\ldots - a_dx_d\\
a_1x^2_1+a_2x^2_2+\ldots+a_{q+1}x^2_{q+1}=\mu_2-a_{q+2}x^2_{q+2}-\ldots - a_dx^2_d\\
..........\\
a_1x^{q}_1+a_2x^{q}_2+\ldots+a_{q+1}x^{q}_{q+1}=\mu_{q}-a_{q+2}x^{q}_{q+2}-\ldots - a_dx^{q}_d\\
\end{array}
\end{equation}
The left hand side of (\ref{eq:Prony.fol.new.2}) is the Vandermonde linear system with respect to $a_1,\ldots,a_{q+1}$.
Hence we can express from (\ref{eq:Prony.fol.new.2}) the amplitudes $a_1,\ldots,a_{q+1}$ via the Cramer rule.
The resulting expressions will be linear in $\mu$ and in $a_{q+2},\ldots,a_d,$ with the coefficients - rational functions in the nodes.
Notice that the denominator is the Vandermonde determinant $V_q(x_1,\ldots,x_{q+1})=\prod_{1\leq i < j \leq q+1}(x_j-x_i)$.

For any fixed $X=(x_1,\ldots,x_d)$ the fiber of $\pi$ over $X$ is an affine subset in ${\cal P}_d$ parametrized by $a_{q+2},\ldots,a_d.$
This completes the proof of Theorem \ref{Thm:Prony.leaves.small.q}. $\square$

Let us stress a special case $q=d-1$. In this case we have

\be\label{eq:Cramer1}
a_j=\frac{1}{V_d(x_1,\ldots,x_d)} \sum_{l=0}^q A_l(X)\mu_l, \ j=1,\ldots,d.
\ee
These expressions remain valid also on the Prony leaves $S_q$ for $q\geq d$.

\subsection{Prony leaves $S_q$ for $q\geq d$}\label{Sec:P.leaves.q.big}

\subsubsection{Projections $S^X_q(\mu)$ of $S_q(\mu)$ onto the nodes subspace}\label{Sec:P.leaves.q.big.x}

Starting with $q=d$ the dimension $2d-q-1$ of the Prony leaves $S_q(\mu)$ is smaller than $d$.
Consequently, the projections $S^X_q(\mu)$ of $S_q(\mu)$ onto the nodes subspace ${\cal P}^X_d$ are proper
subvarieties in ${\cal P}^X_d$. On the other hand, by (\ref{eq:Cramer1}), the amplitudes $a_j$ on $S_q(\mu)$
can be uniquely reconstructed from the nodes $X$ (and from $\mu$). Accordingly, we first describe the equations,
defining the projections $S^X_q(\mu)$ of $S_q(\mu)$ onto the nodes subspace ${\cal P}^X_d$.
To obtain these equations we have to eliminate the amplitudes $a_1,\ldots,a_d$ from the equations (\ref{eq:Prony.fol.new.1}).
This can be achieved by substituting into (\ref{eq:Prony.fol.new.1}) the expressions for $a_j$ from (\ref{eq:Cramer1}).
However, for $d>2$ this leads to rather complicated expressions. Instead we use a modification of the classical solution method of the Prony system. Let

\begin{equation}\label{eq:Prony.fol.new.3}
\begin{array}{c}
\sigma_1(x_1,\ldots,x_d)=-(x_1+\ldots+x_d)\\
\sigma_2(x_1,\ldots,x_d)=x_1x_2+x_1x_3+\ldots+x_{d-1}x_d\\

..........\\
\sigma_d(x_1,\ldots,x_d)=(-1)^d x_1x_2\cdot \ldots \cdot x_{d-1}x_d\\
\end{array}
\end{equation}
be the Vieta elementary symmetric polynomials in $x_1,\ldots,x_d,$. We also put $\sigma_0=1$. Thus $\sigma_j$ are the coefficients of the univariate polynomial

$$
Q(z)=\prod_{j=1}^d(z-x_j)=z^d+\sigma_1z^{d-1}+\ldots+\sigma_d = \sum_{i=0}^d \sigma_{d-i}z^i,
$$
whose roots are the nodes $x_1,\ldots,x_d$.

\smallskip

The following system of $q-d+1$ linear equations for $\sigma_1,\ldots,\sigma_d$ forms a part of the standard
(and classical) linear system for the coefficients of the polynomial $Q$ (see, for instance, \cite{Pro,Pet.Plo,Pot.Tas}. For $q=2d-1$ the complete system is obtained):

\begin{equation}\label{eq:Prony.fol.new.4}
\begin{array}{c}
\mu_{d-1}\sigma_1+\mu_{d-2}\sigma_2+\ldots+\mu_0\sigma_d =-\mu_d\\
\mu_{d}\sigma_1+\mu_{d-1}\sigma_2+\ldots+\mu_1\sigma_d =-\mu_{d+1}\\

..........\\
\mu_{q-1}\sigma_1+\mu_{q-2}\sigma_2+\ldots+\mu_{q-d}\sigma_d =-\mu_{q}\\
\end{array}
\end{equation}
Taking into account that $\sigma_0=1,$ this system can be rewritten as
$$
\sum_{i=0}^d\mu_{l-i}\sigma_i=0, \ l=d,\ldots,q.
$$
System (\ref{eq:Prony.fol.new.4}), being a linear system in variables $\sigma_1,\ldots,\sigma_d$,
forms a nonlinear system in $x_1,\ldots,x_d,$ if we consider $\sigma_j$ as the Vieta elementary symmetric polynomials in
$x_1,\ldots,x_d,$. We denote by $Y_q(\mu) \subset {\cal P}^X_d$ the variety of zeroes of this last system.
\bt\label{thm:Vieta}
For $q\geq d$ the projection $S^X_q(\mu)$ of $S_q(\mu)$ onto the nodes subspace ${\cal P}^X_d$ coincides with $Y_q(\mu)$.
\et
\pr
First we show that for $q\geq d$ system (\ref{eq:Prony.fol.new.1}) implies system (\ref{eq:Prony.fol.new.4}).
Indeed, for each $l=d,\ldots,q$ we obtain, using (\ref{eq:Prony.fol.new.1}), that

$$
\sum_{i=0}^{d}\mu_{l-i}\sigma_i = \sum_{i=0}^{d}\sigma_i\sum_{j=1}^d a_jx^{l-i}_j=\sum_{j=1}^d a_j \sum_{i=0}^{d}\sigma_i x^{l-i}_j=\sum_{j=1}^d a_jx^{l-d}_jQ(x_j)=0,
$$
since each node $x_j$ is a root of $Q(x)$. In other words,
for each $(A,X)\in {\cal P}_d$ satisfying system (\ref{eq:Prony.fol.new.1}),
the component $X$ satisfies system (\ref{eq:Prony.fol.new.4}). We conclude that the projection $S^X_q(\mu)$ of $S_q(\mu)$
onto the nodes subspace ${\cal P}^X_d$ is contained in $Y_q(\mu)$.

\smallskip

To prove the opposite inclusion, let us assume that $X=(x_1,\ldots,x_d)\in Y_q(\mu) \subset {\cal P}^X_d,$ i.e. $X=(x_1,\ldots,x_d)$
satisfies system (\ref{eq:Prony.fol.new.4}). We uniquely define the amplitudes $A =(a_1,\ldots,a_d)$ from the
Vandermonde linear system, formed by the first $d$ equations of system (\ref{eq:Prony.fol.new.1}), according to expressions (\ref{eq:Cramer1}). Now we form a signal
$$
F(x)=\sum_{j=1}^d a_j\delta(x-x_j)=(A, X)\in {\cal P}_d,
$$
which by construction satisfies the first $d$ equations of system (\ref{eq:Prony.fol.new.1}).
It remains to show that the last $q-d+1$ equations of (\ref{eq:Prony.fol.new.1}) are satisfied for $F(x)$.

\smallskip

Consider the rational function $R(z)=\sum_{j=1}^d \frac{a_j}{z-x_j}.$ We have $R(z)=\frac{P(z)}{Q(z)}$ for a certain polynomial $P(z)$ of degree $d-1$ and for
$$
Q(z)=\prod_{j=1}^d(z-x_j)=z^d+\sigma_1z^{d-1}+\ldots+\sigma_d,
$$
where $\sigma_i=\sigma_i(x_1,\ldots,x_d),\ i=1,\ldots,d,$ are, as above, the Vieta elementary symmetric polynomials in $x_1,\ldots,x_d,$.

\smallskip

Developing the elementary fractions in $R(z)$ into geometric progressions, we get

\be\label{eq:Prony.fol.new.5}
R(z)=\sum_{k=0}^\infty \frac{m_k}{z^{k+1}}, \ \ m_k=m_k(F)=\sum_{j=1}^d a_jx_j^k.
\ee
Therefore, the moments $m_k=m_k(F), \ k=0,1,\ldots,$ given by the left hand side $\sum_{j=1}^d a_jx_j^k$ of system (\ref{eq:Prony.fol.new.1}),
are the Taylor coefficients of the rational function $R(z)=\frac{P(z)}{Q(z)},$ with $P(z)$ of degree $d-1$, and $Q(z)$ of degree $d$.
Starting with $k=d$ these Taylor coefficients $m_k$ of $R$ are known to satisfy the recurrence relation

\be\label{eq:Prony.fol.new.6}
m_k=-\sum_{s=1}^d \sigma_s m_{k-l},
\ee
$\sigma_l$ being the coefficients of the denominator $Q(z)$ of $R(z)$.
Since by the choice of the amplitudes $a_j$ the first $d$ equations of system (\ref{eq:Prony.fol.new.1}) are satisfied, we conclude that $m_k=\mu_k, \ k=0,\ldots,d-1.$

\smallskip

Now we use the assumption that system (\ref{eq:Prony.fol.new.4}) is satisfied.
Its equations show that $\mu_k$ satisfy exactly the same recurrence relation till $k=q$.
Since the first $d$ terms are the same, we conclude that in fact $m_k=\mu_k, \ k=0,\ldots,q.$

\smallskip

This means that the entire system (\ref{eq:Prony.fol.new.1}) is satisfied.
Consequently, $F\in S_q(\mu),$ and therefore $X=(x_1,\ldots,x_d)\in S^X_q(\mu) \subset {\cal P}^X_d.$
We conclude that $Y_q(\mu)\subset S^X_q.$ This completes the proof of Theorem \ref{thm:Vieta}. $\square$

\medskip

\noindent {\bf Remark} In the above setting of Theorem \ref{thm:Vieta} we do not make any assumption on the rank of linear system (\ref{eq:Prony.fol.new.4}).
It is easy to give examples of a right-hand side $\mu=(\mu_0,\ldots,\mu_q)$ of (\ref{eq:Prony.fol.new.4})
for which the solutions of this system form an empty set, or an affine subspace $L_q(\mu)$ of any dimension not smaller than $2d-q-1$.
Theorem \ref{thm:Vieta}, as well as Theorem \ref{prop:Hyp.Param} below, remain true in each of this cases.
Compare a detailed discussion of the situation for two nodes $(d=2)$ in Section \ref{Sec:two.nodes} below.

\smallskip

The possible degenerations of system (\ref{eq:Prony.fol.new.4}) are closely related to the conditions of
solvability of Prony system (see, for example, Theorem 3.6 of \cite{Bat.Yom3}), and the discussion thereafter.
Both these questions are very important in the robustness analysis of the Prony inversion, but we do not discuss them here.

\subsubsection{Parametrization of Prony leaves $S_q$ for $q\geq d$}\label{Sec:Param.P.leaves.q.big}

Theorem \ref{thm:Vieta} allows us to construct an explicit parametrization of the Prony leaves $S_q(\mu), \ q\geq d$.
It is enough to produce a parametrization of the projections $S^X_q(\mu)$, since the amplitudes $a_j$ are expressed through the nodes $x_j$
and $\mu$ via formulas (\ref{eq:Cramer1}).

\smallskip

Essentially, we follow the classical solution method of Prony systems,
splitting it into two steps: first, solving a {\it linear system} (\ref{eq:Prony.fol.new.4})
with respect to the variables $\sigma_i$, and then finding the roots of a {\it univariate polynomial} $Q(z)$ with the coefficients $\sigma_i$.

\smallskip

Let $V_d\cong {\mathbb R}^d$ be the space of the coefficients $\sigma=(\sigma_1,\ldots,\sigma_d)$
of the polynomials $Q(z)$ (which we identify with the space of the polynomials $Q$ themselves).
Let $\mu=(\mu_0,\ldots,\mu_q)$ be given. Equations (\ref{eq:Prony.fol.new.4})
define an affine subspace $L_q(\mu)\subset V_d$, which is generically of dimension $2d-q-1$
(but, depending on $\mu$, $L_q(\mu)$ may be empty, or of any dimension not smaller than $2d-q-1$).

\smallskip

Consider a subset $H_d\subset V_d$, consisting of {\it hyperbolic} polynomials $Q$, i.e. of those $Q(z)=z^d+\sigma_1z^{d-1}+\ldots+\sigma_d$
with all the roots real (and pairwise different - so we exclude the boundary).
Hyperbolic polynomials correspond to some of the connected components of the complement in $V_d$
of the discriminant set $\Delta_d\subset V_d$. The set $H_d$ is important in many problems,
and it was intensively studied (see, as a small sample, \cite{Arn,Kos} and references therein).
We denote $L^h_q(\mu)$ the intersection of $L_q(\mu)$ and the set $H_d$ of hyperbolic polynomials.


\bd\label{def:root.vieta.map}
The ``root mapping'' $RM_d:H_d\to {\cal P}^X_d$ is defined by
$$
RM_d(Q)=X=(x_1,\ldots,x_d)\in {\cal P}^X_d,
$$
where $x_1<x_2< \ldots < x_d$ are the ordered roots of the hyperbolic polynomial $Q(z)\in H_d$.

\smallskip

The ``Vieta mapping'' ${\cal V}_d: {\cal P}^X_d\to H_d$ is defined by
$$
{\cal V}_d(x_1,\ldots,x_d)=\left (\sigma_1(x_1,\ldots,x_d),\ldots, \sigma_d(x_1,\ldots,x_d)\right ),
$$
where $\sigma_i=\sigma_i(x_1,\ldots,x_d),\ i=1,\ldots,d,$ are the Vieta elementary symmetric polynomials in $x_1,\ldots,x_d,$.
\ed
Clearly, on $H_d$ the root mapping $RM_d$ is regular, and $RM_d={\cal V}^{-1}_d$.
Therefore both the mappings $RM_d:H_d\to {\cal P}^X_d$ and its inverse ${\cal V}_d: {\cal P}^X_d\to H_d$ provide an isomorphism between $H_d$ and
${\cal P}^X_d$.

\smallskip

Now we have all the tools required to describe the parametrization of the Prony leaves:

\bt\label{prop:Hyp.Param}
The mapping $RM_d:H_d\to {\cal P}^X_d$ transforms isomorphically the affine sets $L^h_q(\mu)$ into the Prony leaves $S^X_q(\mu)$,
while its inverse ${\cal V}_d$ transforms back $S^X_q(\mu)$ into $L^h_q(\mu)$. For each $\mu$ the mapping
$$
RM_d: L^h_q(\mu)\to S^X_q(\mu)
$$
provides an isomorphic parametrization of the Prony leaves $S^X_q(\mu)$.

\smallskip

In other words, the Prony leaves $S^X_q(\mu)\subset {\cal P}^X_d$ are parametrized by the hyperbolic polynomials
$Q\in L^h_d(\mu), \ Q(z)=z^d+\sigma_1z^{d-1}+\ldots+\sigma_d,$ via associating to $Q$ its ordered roots $X=(x_1,\ldots,x_d)\in {\cal P}^X_d$.
\et
\pr
By Theorem \ref{thm:Vieta} the Prony leaf $S^X_q(\mu)\subset {\cal P}^X_d$ consists of all $X=(x_1,\ldots,x_d)\in {\cal P}^X_d$ satisfying equations
(\ref{eq:Prony.fol.new.4}). In other words, $S^X_q(\mu)$ consists of all $X$ for which ${\cal V}_d(X)\in L_q(\mu).$

\smallskip

 On the other hand, since all the nodes of the signals $F=(A,X)\in {\cal P}_d$ are real and pairwise different,
 the necessary and sufficient condition for $\sigma$ to have a form $\sigma=\sigma(X)$ is that all the roots of the polynomial
 $Q(z)=z^d+\sigma_1z^{d-1}+\ldots+\sigma_d$ be real and pairwise different, i.e $Q\in L^h_q(\mu)$.
 As a conclusion, associating to each $Q\in L^h_q(\mu)$ its ordered roots $X=(x_1,\ldots,x_d)\in {\cal P}^x_d$ provides the required parametrization of $S^X_q$. $\square$

\smallskip

An immediate consequence of Theorem \ref{prop:Hyp.Param} is that the Prony leaves $S_q(\mu)$ are smooth algebraic
submanifolds in ${\cal P}_d$. We expect that the results of \cite{Arn,Kos,Kur,Par,Fro.Sha} will be relevant in further investigation of the geometry and topology of the Prony leaves.


\subsection{Prony leaves $S_q$ in the case of two nodes}\label{Sec:two.nodes}

Here we illustrate the results of Sections \ref{Sec:P.leaves.q.small} and \ref{Sec:P.leaves.q.big},
providing a complete description of the Prony leaves in the case of two nodes, i.e. for $d=2$.

\smallskip

For $q=0$ and $\mu=(\mu_0)$ the leaves $S_0(\mu)$ are three-dimensional hyperplanes in ${\cal P}_2\cong {\mathbb R}^4$, defined by the equation $a_1+a_2=\mu_0$.

\smallskip

For $q=1=d-1$ and $\mu=(\mu_0,\mu_1)$ the leaves $S_1(\mu)$ are two-dimensional subvarieties in ${\cal P}_2$, defined by the equations

\be\label{eq:Prony.fol.d.2.1}
a_1+a_2=\mu_0, \ \ a_1x_1+a_2x_2=\mu_1.
\ee
This gives

\be\label{eq:Prony.fol.d.2.3}
a_1=\frac{\mu_0x_2-\mu_1}{x_2-x_1}, \ a_2 = \frac{-\mu_0x_1+\mu_1}{x_2-x_1},
\ee
which is a special case, for $d=2$, of expressions (\ref{eq:Cramer1}).

\medskip

Consider now the case $q=2=2d-2,$ and $\mu=(\mu_0,\mu_1,\mu_2)$. Here the leaves $S_2(\mu)$ are (generically) algebraic curves in ${\cal P}_2$, defined by the equations

\be\label{eq:Prony.fol.d.2.4}
a_1+a_2=\mu_0, \ \ a_1x_1+a_2x_2=\mu_1, \ \ a_1x^2_1+a_2x^2_2=\mu_2,
\ee
For the corresponding curve $S^X_2(\mu)$ in the nodes space ${\cal P}^X_2\cong {\mathbb R}^2$ we obtain from Theorem \ref{thm:Vieta}
the equation $\mu_{1}\sigma_1+\mu_{0}\sigma_2 =-\mu_2,$ or
\be\label{eq:Prony.fol.d.2.5}
\mu_0x_1x_2-\mu_1(x_1+x_2)+\mu_2=0.
\ee
This equation leads to three different possibilities:

\smallskip

\noindent 1. If $\mu_0\ne 0$, then the curve $S^X_2(\mu)$ is a hyperbola

\be\label{eq:Prony.fol.d.2.6}
(x_1-\frac{\mu_1}{\mu_0})(x_2-\frac{\mu_1}{\mu_0})+\frac{\mu_0\mu_2-\mu_1^2}{\mu_0^2}=0,
\ee
which is non-singular for $\mu_0\mu_2-\mu_1^2\ne 0$, and degenerates into two orthogonal coordinate lines,
crossing at the diagonal $\{x_1=x_2\},$ for $\mu_0\mu_2-\mu_1^2=0$.

\smallskip

\noindent 2. If $\mu_0=0$, but $\mu_1\ne 0$ then the curve $S^X_2(\mu)$ is a straight line

\be\label{eq:Prony.fol.d.2.7}
x_1+x_2=\mu_2/\mu_1.
\ee
\noindent 3. Finally, if $\mu_0=\mu_1=0$, but $\mu_2\ne 0$ then the curve $S^X_2(\mu)$ is empty,
and for $\mu_0=\mu_1=\mu_2=0$ it coincides with the entire plane ${\cal P}^x_2$. Compare a discussion in the remark after Theorem \ref{thm:Vieta}.

\medskip

\begin{figure}
	\centering
	\includegraphics[scale=0.6]{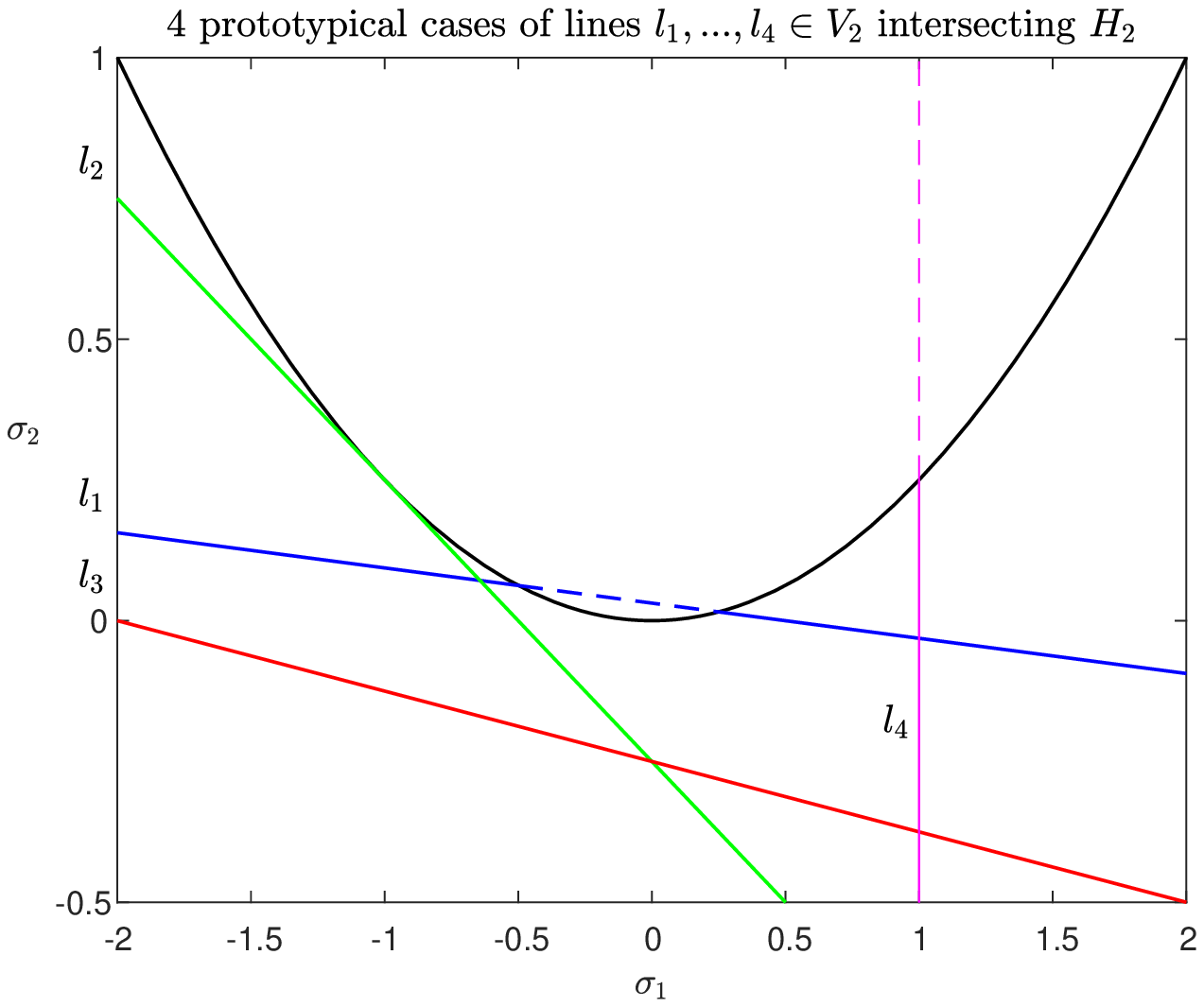}
	\includegraphics[scale=0.75]{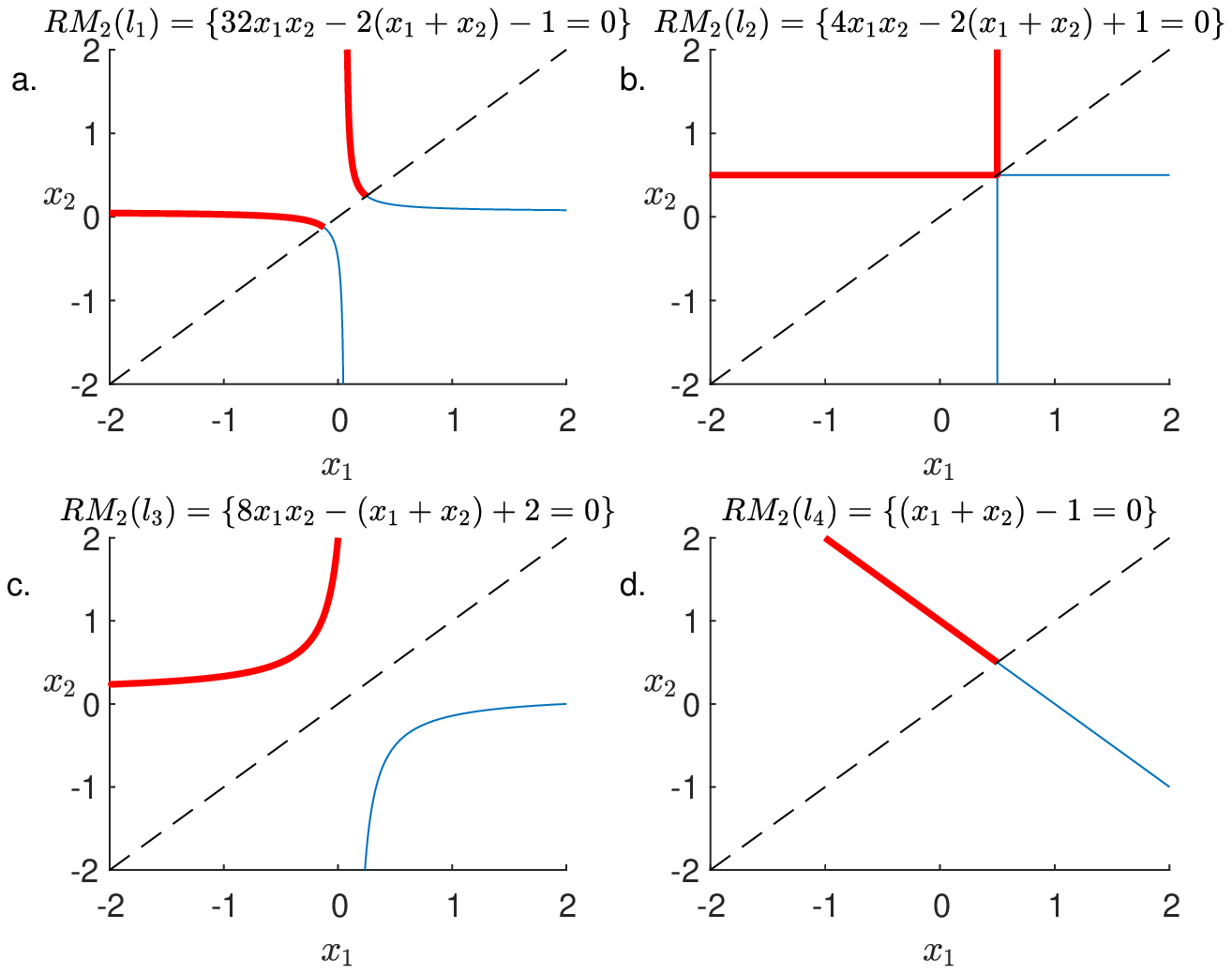}
	\caption{Visualized is the isomorphism $RM_d:H_d\to {\cal P}^X_d$ acting on 4 prototypical lines $l_1,\;l_2\;,l_3\;,l_4 \in V_2$
	intersected with $H_2$  (the open set outside the parabola on the upper figure).
	In the bottom figure, the highlighted parts in the subplots a,b,c,d are the images, under $RM_2$, of $l_1,\;l_2\;,l_3\;,l_4$ intersected with $H_2$, respectively.}
	\label{fig.iso}
\end{figure}

It is instructive to interpret the cases (1-3) above in terms of the relative position,
with respect to the set $H_2$ of hyperbolic polynomials $Q$, of the straight line $L_2(\mu)$.
This line is defined in the space $V_2$ of the polynomials $Q(z)=z^2+\sigma_1 z+\sigma_2$ by system (\ref{eq:Prony.fol.new.4}),
i.e. by the equation $\mu_{1}\sigma_1+\mu_{0}\sigma_2 =-\mu_2$. Figure \ref{fig.iso} illustrates possible positions of the line $L_2(\mu)$ with respect to the set $H_2$ of hyperbolic polynomials.

\smallskip

The discriminant $\Delta(\sigma_1,\sigma_2)=\sigma^2_1-4\sigma_2$ of $Q(z)=z^2+\sigma_1 z+\sigma_2$ is positive for $Q\in H_2$.
Therefore $H_2$ is the part under the parabola $P=\{\sigma_2=\frac{1}{4}\sigma^2_1\}$ in $V_2$. (Compare Figure \ref{fig.iso}).
The case $\mu_0\ne 0$ corresponds to the lines $L_2(\mu)$, nonparallel to the $\sigma_2$-axis of $V_2$.
These lines may cross the parabola $P$ at two points (line $l_1$ on Figure 3), at one point, if tangent to $P$ (line $l_2$ on Figure 3),
or they may not cross $P$ at all, and then they are entirely contained in $H_2$ (line $l_3$ on Figure 3).
These cases correspond to $\mu_0\mu_2-\mu_1^2<0$, $\mu_0\mu_2-\mu_1^2=0$ and $\mu_0\mu_2-\mu_1^2>0$, respectively.

\smallskip

For the line $L_2(\mu)$ crossing the parabola $P$ at two points the corresponding hyperbola $S^X_2(\mu)$
crosses the diagonal in the plane ${\cal P}^X_2$, i.e. it contain a collision of the nodes $x_1,x_2$ (Figure 3, a).

\smallskip

For the line $L_2(\mu)$ tangent to the parabola $P$, the corresponding hyperbola $S^X_2(\mu)$
degenerates into two orthogonal coordinate lines, crossing at a certain point on the diagonal $\{x_1=x_2\},$ (Figure 3, b).

\smallskip

For the line $L_2(\mu)$ entirely contained in $H_2$ the corresponding hyperbola $S^X_2(\mu)$ does not cross the diagonal $\{x_1=x_2\},$
and so it does not lead to the nodes collision (Figure 3, c).

\medskip

For $\mu_0=0$, but $\mu_1\ne 0$, the lines $L_2(\mu)$ are parallel to the $\sigma_2$-axis of $V_2$,
they cross the parabola $P$ at exactly one point (line $l_4$ on Figure 3).
The corresponding curve $S^X_2(\mu)$ is a straight line $x_1+x_2=-\frac{\mu_2}{\mu_1}$ (Figure 3, d).

\smallskip


\section{Some open questions}\label{Sec:Open.Quest}
\setcounter{equation}{0}

The results presented in Section \ref{Sec:Sum.Exs.Results} illustrate the role of the
Prony leaves in the analysis of the error amplification.
The main open problems in the line of this paper concern the structure of the Prony leaves in the areas not covered by the inverse function theorem
(Theorem \ref{thm:coord.moments} above). These areas are collision singularities, on one side,
and ``escape to infinity'' on the other. Both scenarios are frequent in numerical simulations,
but we concentrate on the collision singularities. Let us pose some specific problems in this direction:

\smallskip

\noindent {\it 1. Description of the geometry of the nodes $x_1,\ldots,x_d$ on the Prony leaves $S^x_q(\mu)$ near the collision singularities}.
Presumably, this question can be split into two: investigation of the intersection of the affine varieties $L_q(\mu)\subset V_d$
with the boundary of the hyperbolic set $H_d$, and investigation of the behavior of the root mapping $RM$ near the boundary of $H_d$.

\smallskip

We expect that some classical and more recent results on hyperbolic polynomials,
Vandermonde varieties (and, more generally, on real roots of polynomials and related topics - see \cite{Arn,Fro.Sha,Kos,Kur,Par}) can be relevant.
In particular, in \cite{Kur1} some specific straight lines $l \subset V_d$ are described,
which are entirely contained in $H_d$. Can these lines $l$ appear as the lines $L_{2d-2}(\mu)$ for some $\mu$?

\smallskip

On the other hand, in \cite{Par} smooth selections of real roots in families of polynomials $Q$ are described.
We expect that this description can be relevant in the study of our families $L_q(\mu)$.
``Quantitative'' Lojasievicz-type inequalities may also be useful (see \cite{Kur} and references therein).

\medskip

\noindent {\it 2. Description of the behavior of the amplitudes $a_1,\ldots,a_d$ on the Prony leaves $S^X_q(\mu)$ near the collision singularities}.
We expect that this question can be treated via some methods of the classical Moment theory,
in combination with the techniques of the ``bases of finite differences'' developed in \cite{Bat.Yom3,Yom}.

\smallskip



\medskip

\noindent {\it 3. Extending the description of the Prony leaves, and of the error amplification patterns, to multi-cluster nodes configurations}.
 This is a natural setting in robust inversion of the Prony system.
 In most practical methods separate clusters are first approximated each by a single node,
 thus forming a ``reduced Prony system''. It is important to estimate the accuracy of such an approximation.
\smallskip

Because of the role of the Prony leaves in the analysis of the error amplification patterns,
a natural question is: {\it To what extent the Prony leaves of the reduced Prony system approximate the leaves of the ``true'' multi-cluster system?}


\medskip

\bibliographystyle{amsplain}

\end{document}